\documentclass[leqno,fleqn,a4paper]{amsart}

\usepackage{destemacros}

\begin{document}

\title[Unofficial history of a joint work with Dieter Happel]{Unofficial history of a joint work with Dieter Happel and of two unexpected quotations}


\author{Gabriella D'Este}
\address{Department of Mathematics, University of Milano \\ 
	Via Saldini 50, 
	Milano, Italy
}
\email{gabriella.deste@unimi.it}

\dedicatory{Dedicated to the memory of Dieter Happel}


\date{}

\begin{abstract}
This survey contains a recollection of results, problems
and conversations which go back to the early years
of Representation Theory and Tilting Theory. 
\end{abstract}
\maketitle

\section*{Introduction}
I have many reasons to be very grateful to Dieter Happel.  In this 
note I will describe some facts still vivid in my memory, as a 
personal account of his mathematical vision and (at the same time) 
of his attention to the work of other people, even those he had 
never previously met.  Apart from informal conversations with some 
colleagues in Italy, I remember just one meeting when I said a few 
words about that, namely in Bielefeld (December 2010).  Indeed, in
that occasion, I used the last minutes of my talk to speak about two 
examples of very special and unexpected quotations coming from 
Bielefeld.  The first example is that of an Auslander-Reiten with  28  
vertices, contained in Ringel's survey  \cite[page 93]{R1}.   On the other
hand, the second example, that of a talk, is contained in Happel-Ringel's paper on derived categories  \cite[pages  164  and  180]{HR}.   
This note is organized as follows.  In Section 1, I will recall that I 
met Representation Theory through Happel's handwritten notes  \cite{H}, 
based on a course given by Ringel in Bielefeld.  In Section 2 , I will 
describe the unofficial history of my joint work  \cite{DH}  with Happel.  
Next, in the first part of Section 3 , I will collect some remarks and/or 
conjectures about Auslander-Reiten quivers and sequences.  To this 
end, I will use many pictures concerning either  the quiver with  28  
vertices, mentioned in  \cite{R1}, or other quivers.  Finally, at the end of 
Section 3, I will sum up the prehistory of my Oberwolfach talk, 
mentioned in  \cite{HR}.

\section{Bielefeld, October 1979 :  Happel's notes on
Representation Theory}  
I first met Dieter Happel on October 1979  during my first 
day at Bielefeld Universitaet.  On my second day there
Claus Ringel gave me a copy of the notes of his lectures
of the previous semester written by Happel.  The long 
title of these notes is  ``Vorlesungausarbeitung 
Darstellungstheorie endlich - dimensionaler Algebren''.
I should confess that I read these notes very slowly, 
and with the help of my two dictionaries. Both of them,
a rather small one and a medium sized volume, contained
many useful words, but I was unable to find those typically
German long words, formed by gluing together two or 
three words in some order.  Hence I often had to guess 
possible meanings and/or make conjectures.  It was a new 
experience with respect to both the mathematical subject
and the language.  A direct proof of the last assertion is the
exercise--book, where I wrote the German words I hoped 
to learn.  There are both technical terms from Representation 
Theory and more common words from daily life.  After many 
years, in 1994 , my  10  hours course (``Introduction to the 
representation theory of finite dimensional algebras''), addressed 
to PhD students of the University of Padova,  was nothing else 
but a presentation with slides of a rather small part of  Happel's   
Vorlesungausarbeitung, translated into Italian.  Moreover, I have 
always suggested these notes to colleagues interested to a 
complete introduction to quivers, for both theory and applications, 
and sometimes - the other way round - applications and theory.  
A look at the index indicates the presence of both these aspects
in Happel's notes.  Indeed, the second section (pages  7--19)  is 
dedicated to examples.  Next, the fifth section (pages  39--58)  
is dedicated to the  Auslanderkonstrucktion, that is the direct 
construction of the dual of the transpose  $\tau(M)$   of an 
indecomposable non projective module  $M$, denoted there by  $A(M)$, 
as a tribute to Auslander.  After a section on  $\Ext$, there is a section 
describing the central theorem of the Auslanderkonstruktion and 
some of its applications (Der Hauptsatz und Andwendungen, pages  
64 -73), ending with the pictures of two Auslander-Reiten quivers.  
Finally, the last section deals with the proof of the main theorem  
(Der Beweis des Hauptsatzes, pages 78--88).  The order of the 
various sections and the algorithmic approach mentioned above 
were very important for me at that time, and later on.  Of course, 
I agree that  ``the whole is more than its pieces'', or that  ``the 
Western Wall is more than its pieces'', as I read some years ago 
in the poster of a conference in Analysis in Israel.  However, my 
belief and/or experience is that sometimes also small pieces may 
be useful to understand the whole.  For instance, the maps  
$\tau (M)\To X$   and $X\To M$ of an Auslander--Reiten 
sequence of the form   $0\To\tau(M)\To X\To M\To 0$
are more important, from the functorial point of view, than the 
irreducible maps (between indecomposable modules), the small 
ingredients they are made of, but much more complicated to 
compute and/or guess.  Here I have written ``guess'' for several 
reasons.  First of all, only a few irreducible maps (between 
rather special indecomposable modules) are  well - known 
on the theoretical level.  Next, a kind of concealed ``topology'' 
and geometric symmetries, concerning the shape of close 
irreducible maps already computed, often suggest what should 
be the shape of the still unknown irreducible maps  $\tau(M)\To X(i)$
and  $X(i)\To M$ , where the  $X(i)$'s  are the indecomposable 
summands of  X .  To my astonishment, a mention of  ``intuition''  
shows up in the following final remarks of  Gabriel's paper  \cite[page 66]{G} : ``Since then, various specialists like Bautista, Brenner, 
Butler, Riedtmann\dots  have hoarded a few hundred examples in 
their dossiers, thus getting an intuition which no theoretical 
argument can replace''.

\section{Bielefeld, December 9th , 1988 : A talk at the 
Darstellungstheorie Seminar and what happened 
next.}   
I gave a talk in Bielefeld entitled  ``Some remarks on 
representable equivalences'', which is also the title of
\cite{D} .  Concerning this paper,  I am (and was) very 
grateful to the Editors of the volume ``Topics in Algebras'', 
and in particular to Professor Daniel Simson, for asking 
me to submit a paper, an unexpected opportunity. Indeed, 
the delay of my first flight to Warsaw lead me to cancel 
at the last moment my participation to the conference 
organized by the Banach Center in Warsaw in May 1988, 
at a very difficult time for Poland.   
Back to my Bielefeld seminar, the equivalences I presented 
there are the ones studied by Menini and Orsatti and are 
involved in their Representation Theorem \cite{MO}, a generalization 
of Fuller's Theorem \cite{F}.  The aim of my talk was to show toy 
examples of these equivalences in well - behaved situations, 
by dealing with algebras of finite representation type.  During 
my talk I could immediately answer a question by Happel  
(``Is it extension closed?'')  on one of the classes involved in 
the equivalences studied by Menini and Orsatti in \cite{MO}.  On the 
other hand, I had no answer for other natural questions and/or 
conjectures on the relationship between the  $*$-modules 
considered in  \cite{MO}, and classical tilting modules.  In particular, 
I hoped to find a finite dimensional  $*$-module faithful, but not 
tilting.  The introduction of  \cite{D}  contains the following remark:
``Up to now, we do not know whether or not there exist a finite 
dimensional algebra  A  and a  $*$-module    $_AM$   such that     $_AM$                                                              
is not a  ``disguised'' tilting module, that is     $_{\bar A}M$  is not a tilting
module, where    $\bar A  = A / \ann   _AM$.''   On the other hand, the 
 last words of  \cite{D}  are as follows: 
 ``The proof of  Corollary 6  also  shows that the algebra  $A$  
given by the quiver  $\circ\To\circ\To\circ$   has the following property: 
if  $_AM$ is a multiplicity-free  $*$-module and  $\bar A  = A / \ann   _AM$,
then $_{\bar A} M$    is tilting module.  As already observed in the
Introduction, we do not know any finite dimensional algebra 
without this property''.  
I viewed this last result as a pathological case, and not  as an 
indication of a general property.  On the other hand, Happel 
immediately believed that (over a finite dimensional algebra) 
faithful  $*$-modules and tilting modules always coincide.  I recall 
that he told me more or less the following: ``We should prove this 
fact, and write a joint paper for the journal of Padova University'' 
(the ``Rendiconti del Seminario Matematico dell'Universita' di 
Padova'').   This is the first part of the unofficial history of  \cite{DH}.  
At the same time, it is an example of two decisive factors (somehow 
in the reverse order with respect to the usual way of thinking):
\begin{enumerate}
 \item An overall vision of what should be true.
 \item A winning strategy to prove it, by taking into account 
      all possible contributions of the present and future
      mathematical community.
\end{enumerate}
The second part of this unofficial history is the big role played by 
the copy of the manuscript that Riccardo Colpi gave me shortly 
before my departure.  He told me that it was still work in progress, 
and we had never spoken together in advance about his research.  
I planned to look at Colpi's note only after my visit to Bielefeld.  
However, I putted the envelope in my bag.  When Happel asked 
me whether I had any paper and/or preprint on  $*$-modules of 
Italian colleagues, I said that I had a copy of a preliminary version, 
still in Italian, of the work a PhD student.  After a short look at 
Colpi's handwritten pages, Happel told me that the results 
``looked very interesting''.  Instead of waiting for a final paper 
written in English, he asked me to translate the whole manuscript 
during my short visit to Bielefeld.  To do this, I spent more time 
than expected, because at the same time I wanted to understand 
what I was translating.   As Happel had immediately realized, 
Colpi's notes contained many new ideas and useful tools to deal 
with  $*$-modules and to prove his conjecture on representable 
equivalences.  An evidence of this fact is that we quoted three 
different results of  \cite{C}  
%
(that is,  Corollary 4.2 , Proposition 4.3  
and Theorem 4.1)  in the proof the main theorem in  \cite{DH}.

\section{Bielefeld, December 12th, 2010 :  A colloquium - style 
talk ending with two examples of more than complete 
references (and some speculations about them).}
Thanks to Claudia Koehler, I had the unexpected opportunity to 
present some general ideas on my work to the young audience 
of the Workshop ``Women in Representation Theory: Selfinjective 
Algebras and Beyond'' (December  10--12, 2010), open also 
to  ``senior mathematicians'', as the organizer wrote me.  This was 
the suggestion contained in Claudia's answer to my questions about 
the style of my talk:  ``It would be good if you could give a more 
general talk about something from your research.  If possible 
understandable for  PhD students.''   That's why I prepared a 
personal survey about recent and less results obtained by means 
of the techniques learnt in Bielefeld.  At the end of my talk I read 
the following remarks in the final part of  Ringel's introduction to \cite{R2}:  
``We have tried to give as complete references as possible at the 
end of each chapter, and we apologize for any omission.  Of course, 
it would be easy to trace any omission of a reference to a paper 
which has appeared in print; however we should point out that some 
general ideas which have influenced the results and the methods 
presented here, are not available in official publications, or not even 
written up.''   
I realize that these mathematical and not mathematical reflections 
are very important.  However, the main reason why I mentioned 
these remarks is that I want to recall two unexpected examples 
of quotations (coming from Bielefeld) of quite different results  I 
had obtained (in Bielefeld) long ago.  The first one, concerning an 
Auslander-Reiten quiver, is contained in Ringel's survey  \cite{R1}.  
On the other hand, the second one, concerning more or less large 
modules, is contained in Happel-Ringel's paper  \cite{HR}.  I view 
these two quotations as examples of references which are not only 
complete, but much more than complete.  Both of them are closer 
to true transfigurations (from the real world to an ideal world) 
than to usual references.  Concerning the first case, that of an 
Auslander-Reiten quiver, I will try to give an idea of the big 
gap between real world and ideal world by the means of several 
pictures.  On the other hand, concerning the second case, the 
gap between my result and its presentation is so big that I do not
try to describe any connection between different points of view 
of the same mathematical object.  That's why will copy the short 
account prepared for the Bielefeld Workshop.  I believe that real 
facts speak for themselves.

\subsection{An Auslander-Reiten quiver (elegant form and na•ve     
prehistory).}  The first example mentioned in my Bielefeld talk 
was that the Auslander-Reiten quiver of the algebra, say   $R$, 
given by the quiver 
$\xymatrix{\coltre{}{\bullet}{1} \ar[r]^a & \coltre{}{\bullet}{2} \ar@(ur,dr)[]^{b} }$
	with relations $b^4 = 0$ and $b^2 a = 0$,
%
%
described in Ringel's paper  \cite[page 93]{R1}.  I take from Ringel's 
home page  \cite[Abbildung 3]{R3} the following picture of the 
Auslander-Reiten quiver of  $R$. 
%
\begin{center}
\includegraphics[width=80mm]{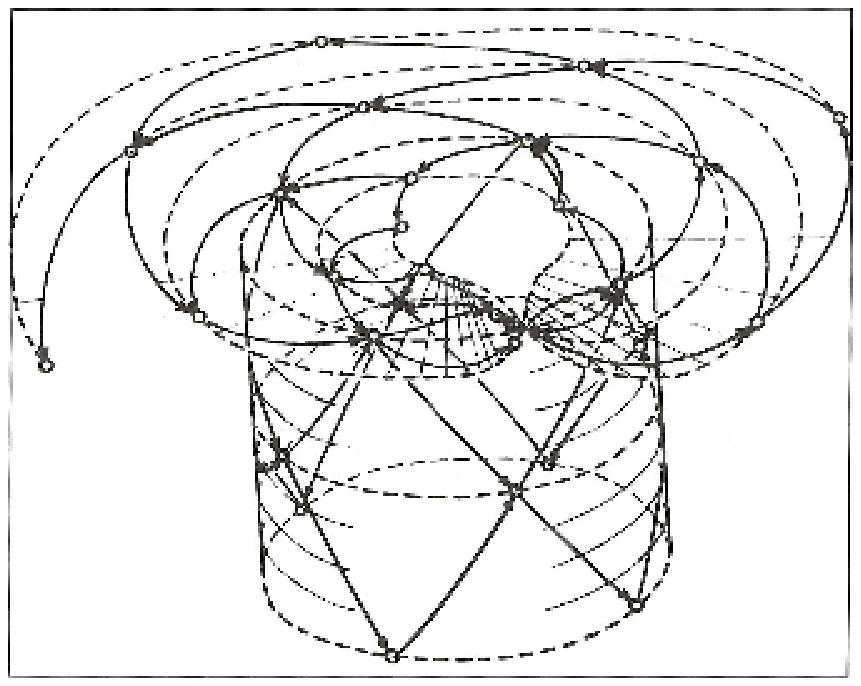}
\end{center}
%
Only after many years -at the end of my talk at a Beijing conference  
\cite{FD}- I expressed for the first time my great surprise of finding a 
quotation (with thanks!) to a quiver I had computed by direct 
calculations and illustrated in a na\"{\i}ve and simple way.  Indeed there 
is an enormous difference between the topological object (formed a 
tube with a M\"obius strip upstair) presented in Ringel's paper and 
home page, and my flat picture of the same Auslander-Reiten quiver 
(formed by two unstable orbits and five stable orbits).  I will need 
several pictures to illustrated some elementary and combinatorial 
aspects of the Auslander-Reiten quiver in my big sheets.  In this 
way I will give a concrete proof of the big gap between the result 
I had and its ideal form presented in  \cite{R1}.  This big gap may be 
an evident example of the ``variety of different appearances'' which 
show up ``in mathematics and in daily life'' , as observed in the 
section ``About Us'' of the home page of  \cite{CRC701}.  
With the usual notation, let  $P(x)$  and  $I(x)$  denote the indecomposable 
projective and injective modules corresponding to the vertex  $x$.  
Moreover, by replacing an indecomposable module which is 
neither projective nor injective by its dimension type, the next 
picture describes the two unstable  $\tau$-orbits of the above 
Auslander-Reiten quiver, containing   3  and  5   indecomposable 
modules respectively, as well as the  6  indecomposable stable 
modules   $X$  in the  $\tau$-orbit of  the radical of  $P(2)$, that is the 
indecomposable module of dimension type  $(0,3)$, such that 
there is an irreducible map of the form  $X\To Y$   or   
$Y\To X$   for some indecomposable unstable module  $Y$.  

\begin{equation*}
	\scalebox{0.55}{
	\xymatrix{
	&&&&&& P(1)\ar[dr] && (1,4)\ar[dr] && I(1) \\
	& P(2)\ar[dr] && (2,2)\ar[dr] && (0,2)\ar[dr]\ar[ur]  && (2,6)\ar[dr]\ar[ur]  && I(2)\ar[dr]\ar[ur]  \\
	(0,3)\ar[ur] && (2,6)\ar[ur] && (2,4)\ar[ur]  && (1,6)\ar[ur] &&(3,6)\ar[ur] &&(1,3)
	}}
\end{equation*}

On the other hand, if we identify the left an right vertical lines in 
order to obtain a M\"obius strip, the next picture describes the 
position of the  20  stable indecomposable modules and their 
dimension types.  
%
More precisely, the $\tau$-orbit of the indecomposable 
module of dimension type  (0,3)  contais  8  modules, 
while the $\tau$-orbits of the indecomposable modules 
of dimension type  $(2,8)$, $(2,7)$  and  $(0,1)$  contain  4
modules respectively.

%
\begin{equation*}
	\hspace*{-12mm}\scalebox{0.7}{
	\xymatrix{
	{}{\ar@[yellow]@{=}@*{[|(10)]}[d]} &&&&&&&& {} {\ar@[yellow]@{=}@*{[|(10)]}[d]} \\
	\fcolorbox{yellow}{white}{(2,4)}\ar[dr] {\ar@[yellow]@{=}@*{[|(10)]}[d]} && (1,6)\ar[dr]  && (3,6)\ar[dr]  && (1,3)\ar[dr]  && \fcolorbox{yellow}{white}{(1,5)} {\ar@[yellow]@{=}@*{[|(10)]}[d]}\\
\fcolorbox{yellow}{white}{(1,5)} \ar[r] {\ar@[yellow]@{=}@*{[|(10)]}[d]}  & \fcolorbox{red}{white}{(3,8)} \ar[r]\ar[dr]\ar[ur]   & (2,3)\ar[r]   & \fcolorbox{red}{white}{(2,6)} \ar[dr]\ar[r]\ar[ur]  & (0,3)\ar[r]  & \fcolorbox{red}{white}{(2,5)} \ar[r] \ar[dr]\ar[ur]  &(2,6)\ar[r]  & \fcolorbox{red}{white}{(2,8)} \ar[r] \ar[dr]\ar[ur]  &  \fcolorbox{yellow}{white}{(2,4)}{\ar@[yellow]@{=}@*{[|(10)]}[d]} \\
\fcolorbox{yellow}{white}{\fcolorbox{green}{white}{(2,7)}}\ar[dr] \ar[ur] {\ar@[yellow]@{=}@*{[|(10)]}[d]}  && \fcolorbox{green}{white}{(2,5)} \ar[dr] \ar[ur]  &&  \fcolorbox{green}{white}{(1,2)} \ar[dr] \ar[ur]  && \fcolorbox{green}{white}{(1,4)} \ar[dr] \ar[ur]  && {\fcolorbox{yellow}{white}{\fcolorbox{green}{white}{(2,7)}} }{\ar@[yellow]@{=}@*{[|(10)]}[d]}   \\
& \fcolorbox{blue}{white}{(1,4)} \ar[ur] && \fcolorbox{blue}{white}{(1,1)} \ar[ur] && \fcolorbox{blue}{white}{(0,1)} \ar[ur] && \fcolorbox{blue}{white}{(1,3)} \ar[ur]  &                 
	}}
\end{equation*}

\subsection{ An Auslander-Reiten sequence of the previous Auslander-Reiten quiver (= example where all irreducible maps are 
obvious cancellations or additions).}   By comparing dimension 
types, we deduce from the above pictures that the dimension types 
of the stable modules considered above --that is  (0,1), (0,3), (2,7)  
and  (2,8)-- are the dimension types of exactly one indecomposable 
$R$-module.  This is obvious for the first two cases, which are the 
dimension types of the simple top of   $P(2)$  and of its unique maximal 
submodule.   Let  $L$  and  $M$  denote the above indecomposable modules 
of dimension type  (2,7)  and  (2,8)  respectively.  Then the matrices   
$A$  and  $B$  (resp.  $A'$  and  $B'$), describing the injective map  $a$   from 
$K^2$ to $K^7$ (resp. $K^8$)   and the endomorphism  $b$  of  $K^7$    (resp.  $K^8$   ),
are of the following form:

\begin{equation*}
	A=
	\begin{pmatrix}
		0 & 0 \\
	 	0 & 0 \\
		1 & 0 \\
		0 & 0 \\
		0 & 0 \\
		1 & 0 \\
		0 & 1 
	\end{pmatrix},\qquad 
	B= 
	\begin{pmatrix} 
		0 & 0 & 0 & 0 & 0 & 0 & 0  \\
		1 & 0 & 0 & 0 & 0 & 0 & 0  \\
		0 & 1 & 0 & 0 & 0 & 0 & 0  \\
		0 & 0 & 1 & 0 & 0 & 0 & 0  \\
		0 & 0 & 0 & 0 & 0 & 0 & 0  \\
		0 & 0 & 0 & 0 & 1 & 0 & 0  \\
		0 & 0 & 0 & 0 & 0 & 1 & 0  
	\end{pmatrix}
\end{equation*}

\begin{equation*}
	A'=
	\begin{pmatrix}
	 	1 & 0 \\
		0 & 0 \\
		0 & 0 \\
		1 & 0 \\
		0 & 0 \\
		0 & 0 \\
		1 & 0 \\
		0 & 1 
	\end{pmatrix},\qquad 
	B'= 
	\begin{pmatrix} 
		0 & 0 & 0 & 0 & 0 & 0 & 0 & 0 \\
		0 & 0 & 0 & 0 & 0 & 0 & 0 & 0 \\
		0 & 1 & 0 & 0 & 0 & 0 & 0 & 0 \\
		0 & 0 & 1 & 0 & 0 & 0 & 0 & 0 \\
		0 & 0 & 0 & 1 & 0 & 0 & 0 & 0 \\
		0 & 0 & 0 & 0 & 0 & 0 & 0 & 0 \\
		0 & 0 & 0 & 0 & 0 & 1 & 0 & 0 \\
		0 & 0 & 0 & 0 & 0 & 0 & 1 & 0 
	\end{pmatrix}
\end{equation*}
I wrote these big matrices just to explain with two reasonably small 
pictures the conventions, suggested me by Ringel during my stay in 
Bielefeld, used to describe in a compact (and efficient !) way the
indecomposable representations of the form  
\begin{equation*}
 \xymatrix{\scalebox{0.85}{V(1)} \ar[r]^a & \scalebox{0.85}{V(2)} \ar@(ur,dr)[]^{b} }
\end{equation*}
with  $V(2)\neq 0$, that is different from the simple injective module 
$I(1)$, for the various algebras given by the quiver
\begin{equation*}
 \xymatrix{\coltre{}{\bullet}{1} \ar[r]^a & \coltre{a}{\bullet}{2} \ar@(ur,dr)[]^{b} }
\end{equation*}
 with relations   $b^n = 0$   and   $b^2 a = 0$ for some  $n > 1$.  
 Following Ringel's hints, every square in the next 
pictures indicates a fixed element of a fixed basis of the vector space  
$V(2)$, while the various blocks indicate the uniserial summands of the 
submodule   
\begin{equation*}
 \xymatrix{\scalebox{0.85}{0} \ar[r] & \scalebox{0.85}{V(2)}  \ar@(ur,dr)[]^{b} }.
\end{equation*}
On the other, every black
square indicates a vector   $v$  of the fixed basis belonging to the image
of the injective linear map  $a$.  Finally, every segment connecting finitely 
many small squares, corresponding to vectors   $v(1), v(2),\dots, v(m)$, 
indicates that  $v(1) + v(2) +\dots+ v(m)$  belongs to the image of  $a$.   
Following these conventions, the previous modules  $L$  and  $M$  , of 
dimension type  (2,7)   and  (2,8)  look like as follows.
\subsection*{Picture 1}
\begin{center}
  \begin{tabular}{@{} ccc @{}}
    $
\xymatrix @R=0pc @C=1pc @H=1pc @W=2pc  
	{
	*++[F-]{\phantom{\bullet}}  & {\phantom{\bullet}}\\
	*++[F-]{\phantom{\bullet}} 	& *++[F-]{\phantom{\bullet}}  \\
	*++[F-]{\bullet}  \ar@{-}[r]  & *++[F-]{\bullet}  \\
	*++[F-]{\phantom{\bullet}}  & *++[F*]{\phantom{\bullet}}\\
}$ & {\xymatrix{\\ \txt{,}\\ }} & 
$\xymatrix @R=0pc @C=1pc @H=1pc @W=2pc  
{	&  *++[F-]{\phantom{\bullet}} & \\
	&  *++[F-]{\phantom{\bullet}} &  *++[F-]{\phantom{\bullet}} \\
 *++[F-]{\bullet} \ar@{-}[r]	&  *++[F-]{\bullet} \ar@{-}[r] &  *++[F-]{\bullet} \\
	&  *++[F-]{\phantom{\bullet}} &  *++[F*]{\phantom{\bullet}} \\
}$
\\ 
    $L=(2,7)$ &  & $M=(2,8)$
  \end{tabular}
\end{center}
Consequently, the arrow from  (2,8)  to  (2,7)  indicates
a kind of left cancellation, that is the irreducible 
epimorphism whose simple kernel is the unique simple
summand of   
\begin{equation*}
 \xymatrix{\scalebox{0.85}{0} \ar[r] & {\scalebox{0.85}{$K^8$}} \ar@(ur,dr)[]^{b} }.
\end{equation*}
On the other 
hand, the arrow from  (2,8)  to   (2,4)  (resp.   (1,5))  
describes a kind of central (resp.  right)  cancellation of the     
unique summands of dimension 4 (resp. 3) of the above
 submodule of dimension type (0,8).
Indeed, these
two modules have the following shape:

\subsection*{Picture 2}
\begin{center}
  \begin{tabular}{@{} ccc @{}}
    $
\xymatrix @R=0pc @C=1pc @H=1pc @W=2pc  
	{
	  							& *++[F-]{\phantom{\bullet}}\\
	*++[F-]{\bullet}  \ar@{-}[r]  & *++[F-]{\bullet}  \\
	  & *++[F*]{\phantom{\bullet}}\\
}
$ & {\xymatrix{\\ \txt{,}\\ }} & $\xymatrix @R=0pc @C=1pc @H=1pc @W=2pc  
	{
	  							& *++[F-]{\phantom{\bullet}}\\
	  							& *++[F-]{\phantom{\bullet}}\\
	*++[F-]{\bullet}  \ar@{-}[r]  & *++[F-]{\bullet}  \\
	  & *++[F-]{\phantom{\bullet}}\\
}
$ \\ 
    (2,4) &  & (1,5)
  \end{tabular}
\end{center}

Dually, keeping the above conventions, the indecomposable 
module of dimension type  (3,8) has the following shape

\subsection*{Picture 3}
$$
	\xymatrix @R=0pc @C=1pc @H=1pc @W=2pc  
	{
	  							& *++[F-]{\phantom{\bullet \ \bullet}} & \\
	  							& *++[F-]{\phantom{\bullet \ \bullet}} & *++[F-]{\phantom{\bullet}}\\
	*++[F-]{\bullet}  \ar@{-}[r]  & *++[F-]{\bullet \ \bullet}  \ar@{-}[r]  & *++[F-]{\bullet}  \\
	  & *++[F-]{\phantom{\bullet \ \bullet}} & *++[F*]{\phantom{\bullet}}
	  }
$$

Consequently, the three irreducible maps arriving at (3,8) look like
additions of different types, that is either left or right or
central additions.

\subsection{Comparing (obvious and less obvious) Auslander - 
Reiten sequences.}   We list some other properties of the 
Auslander-Reiten sequence ending at the above
indecomposable module of dimension type  (3,8), briefly 
denoted by
\begin{equation}
\xymatrix@1{0 \ar[r] & (2,8) \ar[r]^f  &  X \ar[r]^g  &  (3,8) \ar[r]  &  0}
\tag{$\sharp$}
\end{equation}
First of all, for every indecomposable non projective module  M  
of dimension type different from  (3,8), we have $\dim_K M + \dim_k\tau(M) < 21$. 
  Hence it is reasonable to expect
that  ($\sharp$)  is one of the most complicated Auslander-Reiten 
sequences of the above Auslander-Reiten quiver.  Secondly,  
$X$  is the direct sum of  3  indecomposable modules, neither 
projective nor injective.  Therefore, by the four-in-the-middle theorem \cite{BB},  $X$  has the largest possible number of 
indecomposable and stable direct summands.  However, even 
in this complicated case, the maps  $f$  and  $g$  in  ($\sharp$)  consist 
of  three irreducible maps which are very natural ``cancellations'' 
and ``additions'' respectively.  Hence, we may roughly speaking 
say that the functorial and global properties of  $f$  and $g$  are 
somehow concealed.  On the other hand the local properties of 
the irreducible components of  $f$  and  $g$  are evident, at least 
``intuitively'', to repeat the words used in Gabriel's remark  \cite[page  66]{G} (see the end of section 1).  In the above sequence  
($\sharp$)  we may see at a glance that all the three irreducible 
components of  $f$  and  $g$  are the obvious ones, by looking at 
exactly one visualization of the first and last non--zero modules.  
This looks like a nice situation, but I do not view it as an 
exceptional one.  Indeed, if  $Z$  is an indecomposable summand 
of the middle term of an Auslander-Reiten sequence, say  $(*)$,  
and  $h$  is an irreducible component of one of the two non--zero 
maps in  $(*)$  arriving at  $Z$  or ending in  $Z$ , then only one of the 
following cases seems to occur.


\begin{description}
\item[The obvious good case] The map  $h$  looks like an evident 
injective or surjective map, that is a kind of translation of the 
unique well-known irreducible maps at the theoretical level,
of the form   $X\To P$ and $I\To Y$, where the modules
$X, P, I, Y$  are indecomposable,  $P$  is projective and  $I$  is
injective. 
\item[The concealed good case] Up to a change of the basis of 
the first or last non-zero terms of  $(*)$, $h$  becomes an obvious
embedding or epimorphism.  Moreover, after that change,  $h$
looks like as a shift of other irreducible maps close to it.
\end{description}

\subsection{A more complicated Auslander-Reiten sequence  
(= example where all irreducible maps become obvious 
cancellations or additions after changing a basis).}
After an example of the good case, given by the previous 
sequence  ($\sharp$),  I will describe an example of a concealed good 
case by means of an Auslander-Reiten sequence, say  ($\sharp\sharp$), 
ending at an indecomposable module of dimension type  (6,24), 
briefly denoted by
\begin{equation}
\tag{$\sharp \sharp$} \xymatrix@1{ 0 \ar[r] & (7,27) \ar[r] & X \ar[r] & (6,24) \ar[r] & 0 }.
\end{equation}

In this case the algebra into the game is the path algebra
given by the quiver   
\begin{equation*}
\xymatrix@1{	\coltre{}{\bullet}{1} \ar[r]^a & \coltre{}{\bullet}{2} \ar@(ur,dr)[]^{b}   } 
\qquad \text{with relations }\quad b^6 = 0\ \text{ and } \ b^2 a = 0.
\end{equation*}

Also in this case the middle term  $X$   is the direct sum of three 
indecomposable non projective summands of dimension types  
(5 , 18),  (3 , 13)   and  (5 , 20)  with the following shapes
respectively:
\subsection*{(5,18)}
$$
	\xymatrix @R=0pc @C=1pc @H=1pc @W=2pc  
	{
	  						&	& *++[F-]{\phantom{\bullet}} & & \\
	*++[F-]{\phantom{\bullet}} 	&	& *++[F-]{\phantom{\bullet}} & *++[F-]{\phantom{\bullet}} & \\
	*++[F-]{\phantom{\bullet}} 	&	& *++[F-]{\phantom{\bullet}} & *++[F-]{\phantom{\bullet}} & *++[F-]{\phantom{\bullet}}\\
	*++[F-]{\bullet} \ar@{-}[r]	& *++[F-]{\bullet} & *++[F-]{\phantom{\bullet}} & *++[F-]{\bullet} \ar@{-}[r] & *++[F-]{\bullet}\\
	*++[F-]{\bullet} \ar@{-}[rr]	&	& *++[F-]{\bullet} \ar@{-}[r] & *++[F-]{\bullet} & *++[F*]{\phantom{\bullet}}\\
	  & & *++[F*]{\phantom{\bullet}} & &
	  }
$$

\subsection*{(3,13)}
$$
	\xymatrix @R=0pc @C=1pc @H=1pc @W=2pc  
	{
	 *++[F-]{\phantom{\bullet}} & & \\
	 *++[F-]{\phantom{\bullet}} & *++[F-]{\phantom{\bullet}} & \\
	 *++[F-]{\phantom{\bullet}} & *++[F-]{\phantom{\bullet}} & *++[F-]{\phantom{\bullet}} \\
	 *++[F-]{\phantom{\bullet}} & *++[F-]{\bullet}\ar@{-}[r] & *++[F-]{\bullet} \\
	 *++[F-]{\bullet} \ar@{-}[r]& *++[F-]{\bullet} & *++[F*]{\phantom{\bullet}} \\
	 *++[F-]{\phantom{\bullet}} & & \\
	  }
$$
\subsection*{(5,20)}
$$
	\xymatrix @R=0pc @C=1pc @H=1pc @W=2pc  
	{
	 *++[F-]{\phantom{\bullet}} & &  & *++[F-]{\phantom{\bullet}}  & \\
	 *++[F-]{\phantom{\bullet}} & *++[F-]{\phantom{\bullet\ \bullet}} &  & *++[F-]{\phantom{\bullet}} & \\
	 *++[F-]{\phantom{\bullet}} & *++[F-]{\phantom{\bullet\ \bullet}} &   
	 & *++[F-]{\phantom{\bullet}} &  *++[F-]{\phantom{\bullet}} \\
	 *++[F-]{\phantom{\bullet}} & *++[F-]{\phantom{\bullet}\ \bullet} \ar@{-}[r]& *++[F-]{\bullet}  
	 & *++[F-]{\phantom{\bullet}} &  *++[F*]{\phantom{\bullet}} \\
	 *++[F-]{\bullet} \ar@{-}[r]& *++[F-]{\bullet\ \bullet} \ar@{-}[rr] & 
	 & *++[F-]{\bullet} \ar@{-}[r] &  *++[F-]{\bullet} \\
	 *++[F-]{\phantom{\bullet}} & &    
	 & *++[F*]{\phantom{\bullet}} &   \\
	  }
$$

We illustrate in the sequel the original shapes (indexed by  $A$ )  of the
first and the last non-zero term of   ($\sharp \sharp$) , and those later obtained 
after a change of the basis (indexed by  $B$  and  $C$ ).

\subsection*{(6,24) A}
$$
	\xymatrix @R=0pc @C=1pc @H=1pc @W=2pc  
	{
	 *++[F-]{\phantom{\bullet}} & &  & *++[F-]{\phantom{\bullet}}  & & \\
	 *++[F-]{\phantom{\bullet}} & *++[F-]{\phantom{\bullet\ \bullet}} &  & *++[F-]{\phantom{\bullet}} & *++[F-]{\phantom{\bullet}}  & \\
	 *++[F-]{\phantom{\bullet}} & *++[F-]{\phantom{\bullet\ \bullet}} &   
	 & *++[F-]{\phantom{\bullet}} &  *++[F-]{\phantom{\bullet}} &  *++[F-]{\phantom{\bullet}} \\
	 *++[F-]{\phantom{\bullet}} & *++[F-]{\phantom{\bullet}\ \bullet} \ar@{-}[r]& *++[F-]{\bullet}  
	 & *++[F-]{\phantom{\bullet}} &  *++[F-]{\bullet} \ar@{-}[r] &  *++[F-]{\bullet} \\
	 *++[F-]{\bullet} \ar@{-}[r]& *++[F-]{\bullet\ \bullet} \ar@{-}[rr] & 
	 & *++[F-]{\bullet} \ar@{-}[r] &  *++[F-]{\bullet} &  *++[F*]{\phantom{\bullet}} \\
	 *++[F-]{\phantom{\bullet}} & &    
	 & *++[F*]{\phantom{\bullet}} &  & \\
	  }
$$

\subsection*{(6,24) B}
$$
	\xymatrix @R=0pc @C=1pc @H=1pc @W=2pc  
	{
	 & & *++[F-]{\phantom{\bullet}}  & *++[F-]{\phantom{\bullet}}  & & \\
	& *++[F-]{\phantom{\bullet}} & *++[F-]{\phantom{\bullet}}  & *++[F-]{\phantom{\bullet}} & *++[F-]{\phantom{\bullet}}  & \\
	  & *++[F-]{\phantom{\bullet}} &   *++[F-]{\phantom{\bullet}}
	 & *++[F-]{\phantom{\bullet}} &  *++[F-]{\phantom{\bullet}} &  *++[F-]{\phantom{\bullet}} \\
	 *++[F-]{\bullet} \ar@{-}[r]& *++[F-]{\bullet} & *++[F-]{\phantom{\bullet}}  
	 & *++[F-]{\phantom{\bullet}} &  *++[F-]{\bullet} \ar@{-}[r] &  *++[F-]{\bullet} \\
	 & *++[F-]{\bullet} \ar@{-}[r] & *++[F-]{\bullet}
	 & *++[F-]{\bullet} \ar@{-}[r] &  *++[F-]{\bullet} &  *++[F*]{\phantom{\bullet}} \\
	 & &   *++[F-]{\bullet} \ar@{-}[r]
	 & *++[F-]{\bullet} &  & \\
	  }
$$
\subsection*{(6,24) C}
$$
	\xymatrix @R=0pc @C=1pc @H=1pc @W=2pc  
	{
	 *++[F-]{\phantom{\bullet}} & &   & *++[F-]{\phantom{\bullet}}  & & \\ 
	 *++[F-]{\phantom{\bullet}} & *++[F-]{\phantom{\bullet \ \bullet}} &   & *++[F-]{\phantom{\bullet}}  & & *++[F-]{\phantom{\bullet}} \\ 
	 *++[F-]{\phantom{\bullet}} & *++[F-]{\phantom{\bullet \ \bullet}} &   & *++[F-]{\phantom{\bullet}}  & *++[F-]{\phantom{\bullet \ \bullet}} & *++[F-]{\phantom{\bullet}} \\ 
	 *++[F-]{\phantom{\bullet}} & *++[F-]{\phantom{\bullet} \ \bullet} \ar@{-}[r]& *++[F-]{\bullet}   & *++[F-]{\phantom{\bullet}}  & *++[F*]{\phantom{\bullet \ \bullet}} & *++[F-]{\phantom{\bullet}} \\ 
	 *++[F-]{\bullet}\ar@{-}[r] & *++[F-]{\bullet \ \bullet} \ar@{-}[rr]&    & *++[F-]{\bullet}  \ar@{-}[r] & *++[F-]{\bullet \ \bullet} \ar@{-}[r] & *++[F-]{\bullet}  \\ 
	 *++[F-]{\phantom{\bullet}} & &   & *++[F*]{\phantom{\bullet}}  & &  
  }
$$

\subsection*{(7,27) A}
$$
	\xymatrix @R=0pc @C=1pc @H=1pc @W=2pc  
	{
	& & *++[F-]{\phantom{\bullet}} & &   & *++[F-]{\phantom{\bullet}} & \\ 
	& *++[F-]{\phantom{\bullet\ \bullet}} & *++[F-]{\phantom{\bullet}} & *++[F-]{\phantom{\bullet\ \bullet}} &   & *++[F-]{\phantom{\bullet}} & \\ 
*++[F-]{\phantom{\bullet}} 	& *++[F-]{\phantom{\bullet\ \bullet}} & *++[F-]{\phantom{\bullet}} & *++[F-]{\phantom{\bullet\ \bullet}} &    & *++[F-]{\phantom{\bullet}} & *++[F-]{\phantom{\bullet}}  \\ 
*++[F*]{\phantom{\bullet}} 	& *++[F-]{\phantom{\bullet\ \bullet}} & *++[F-]{\phantom{\bullet}} & *++[F-]{\phantom{\bullet}\ \bullet} \ar@{-}[r] & *++[F-]{\bullet}   & *++[F-]{\phantom{\bullet}} & *++[F*]{\phantom{\bullet}}  \\ 
*++[F-]{\bullet} \ar@{-}[r]	& *++[F-]{\bullet\ \bullet} \ar@{-}[r] & *++[F-]{\bullet} \ar@{-}[r] & *++[F-]{\bullet\ \bullet}\ar@{-}[rr] &     & *++[F-]{\bullet}\ar@{-}[r] & *++[F-]{\bullet}  \\ 
	& & *++[F-]{\phantom{\bullet}} & &   & *++[F*]{\phantom{\bullet}} &  
  }
$$

\subsection*{(7,27) B}
$$
	\xymatrix @R=0pc @C=1pc @H=1pc @W=2pc  
	{
	*++[F-]{\phantom{\bullet}} & &  & &   & *++[F-]{\phantom{\bullet}} & \\ 
	*++[F-]{\phantom{\bullet}} & *++[F-]{\phantom{\bullet}} & *++[F-]{\phantom{\bullet}} & *++[F-]{\phantom{\bullet\ \bullet}} &   & *++[F-]{\phantom{\bullet}} & \\ 
	*++[F-]{\phantom{\bullet}} & *++[F-]{\phantom{\bullet}} & *++[F-]{\bullet} & *++[F-]{\phantom{\bullet\ \bullet}} &   & *++[F-]{\phantom{\bullet}} & *++[F-]{\phantom{\bullet}} \\ 
	*++[F-]{\phantom{\bullet}} & *++[F-]{\bullet} \ar@{-}[ru] & *++[F*]{\phantom{\bullet}} & *++[F-]{\phantom{\bullet} \ \bullet} \ar@{-}[r]&  *++[F-]{\bullet} & *++[F-]{\phantom{\bullet}} & *++[F*]{\phantom{\bullet}} \\ 
	*++[F-]{\bullet} \ar@{-}[r] & *++[F-]{\bullet} \ar@{-}[rr] &  & *++[F-]{\bullet\ \bullet} \ar@{-}[rr]&    & *++[F-]{\bullet} \ar@{-}[r] & *++[F-]{\bullet} \\ 
	*++[F-]{\phantom{\bullet}} & &  & &   & *++[F*]{\phantom{\bullet}} & \\ 
  }
$$

\subsection*{(7,27) C}
$$
	\xymatrix @R=0pc @C=1pc @H=1pc @W=2pc  
	{
	& *++[F-]{\phantom{\bullet \ \bullet}} & &   & *++[F-]{\phantom{\bullet}} &  &  \\ 
	& *++[F-]{\phantom{\bullet \ \bullet}} & *++[F-]{\phantom{\bullet}} &   & *++[F-]{\phantom{\bullet}} & *++[F-]{\phantom{\bullet \ \bullet}} &  \\ 
	*++[F-]{\phantom{\bullet}} & *++[F-]{\phantom{\bullet \ \bullet}} & *++[F-]{\phantom{\bullet}} &   & *++[F-]{\phantom{\bullet}} & *++[F-]{\phantom{\bullet \ \bullet}} &  *++[F-]{\phantom{\bullet}} \\ 
	*++[F*]{\phantom{\bullet}} & *++[F-]{\phantom{\bullet \ \bullet}} & *++[F-]{\bullet} \ar@{-}[r] &  *++[F-]{\bullet} & *++[F-]{\phantom{\bullet}} & *++[F-]{\phantom{\bullet \ \bullet}} &  *++[F*]{\phantom{\bullet}} \\ 
	*++[F-]{\bullet} \ar@{-}[r] & *++[F-]{\bullet \ \bullet} \ar@{-}[r]& *++[F-]{\bullet}\ar@{-}[rr] &  & *++[F-]{\bullet}\ar@{-}[r] & *++[F-]{\bullet \ \bullet}\ar@{-}[r] &  *++[F-]{\bullet} \\ 
	& *++[F-]{\phantom{\bullet \ \bullet}} & &   & *++[F*]{\phantom{\bullet}} &  &  \\ 
  }
$$


The calculation of   ($\sharp\sharp$)  goes back to the first problem suggested by  
Ringel (with many useful hints and suggestions) at the end of my first 
year in Bielefeld.  By looking at my old pictures of  $(6,24)$  and  $(7,27)$, 
that is  $(6,24) A$  and  $(7,27) A$), only two out of six irreducible maps, 
that is   $(5,18) \To  (6,24) A$   and  $(7,27) A \To (5,20)$, look like 
obvious ``additions'' and ``cancellations''.  In order to see that the same 
holds for the other four irreducible maps, it suffices to use the 
visualizations of type  $B$  and  $C$  of the first and the last non zero 
modules in  ($\sharp\sharp$).
We illustrate in the sequel the behaviour of the three 
reducible maps obtained by making use of the obvious form of the six 
irreducible components of the two non -zero maps in  ($\sharp\sharp$) .  First of all, 
the composition   $(7,27) C  \To (5,18)\To (6,24) A$    acts as 
a left cancellation followed by a left addition.
$$
	\xy
	\xymatrix"A" @R=0pc @C=1pc @H=1pc @W=2pc  
	{
	& *++[F-]{\phantom{\bullet \ \bullet}}  & \\ 
	& *++[F-]{\phantom{\bullet \ \bullet}}  & \\ 
	*++[F-]{\phantom{\bullet}} & *++[F-]{\phantom{\bullet \ \bullet}} &  \\ 
	*++[F*]{\phantom{\bullet}} & *++[F-]{\phantom{\bullet \ \bullet}}  & \\ 
	*++[F-]{\bullet} \ar@{-}[r] & *++[F-]{\bullet \ \bullet}  & 
	 \\ 
	 & *++[F-]{\phantom{\bullet \ \bullet}}  & \\ 
	}
	\POS+(30,-13)
	\xymatrix{
	 &  \\
	 &  \\
	 {\bullet} \ar@{-}["A"ddr] & (5,18) \\
	 & \\
	}
	\POS*\frm{--}
	\endxy
	\qquad \xymatrix{\\ \\ \scalebox{3.}{$\rightarrow$}} \qquad 
	\xy
	\xymatrix"B" @R=0pc @C=1pc @H=1pc @W=2pc  
	{
	 *++[F-]{\phantom{\bullet}}  \\ 
	 *++[F-]{\phantom{\bullet}}  \\ 
	 *++[F-]{\phantom{\bullet}}  \\ 
	 *++[F-]{\phantom{\bullet}}  \\ 
	 *++[F-]{\bullet} \\ 
	 *++[F-]{\phantom{\bullet}}  \\ 
	}
	\POS+(13,-10)
	\xymatrix{
	 &  \\
	 &  \\
	 {\bullet} \ar@{-}["B"dd] & (5,18) \\
	 & \\
	}
	\POS*\frm{--}
	\endxy
$$
Secondly, the composition  $(7,27) B \To (3,13) \To  (6,24) B$
acts as a right cancellation followed by a left addition. 
$$
	\xy
	\xymatrix"B"{
	  (3,13)  \\
	 {\bullet}  \\
	}
	\POS*\frm{--}
	\POS+(13,12)
	\xymatrix"C" @R=0pc @C=1pc @H=1pc @W=0.7pc  
	{
	 &  & *++[F-]{\phantom{\bullet}} & \\ 
	*++[F-]{\phantom{\bullet\ \bullet}}  &  & *++[F-]{\phantom{\bullet}} & \\ 
	*++[F-]{\phantom{\bullet\ \bullet}}  &  & *++[F-]{\phantom{\bullet}} & *++[F-]{\phantom{\bullet}}\\ 
	*++[F-]{\phantom{\bullet}\ \bullet} \ar@{-}[r] & *++[F-]{\bullet} & *++[F-]{\phantom{\bullet}} & *++[F*]{\phantom{\bullet}}\\ 
	*++[F-]{\bullet\ \bullet} \ar@{-}[rr] \ar@{-}["B"uuu] &  & *++[F-]{\bullet} \ar@{-}[r] & *++[F-]{\bullet}\\ 
	& & *++[F*]{\phantom{\bullet}} &  
	}
	\endxy
	\ \xymatrix{\\  \scalebox{3.}{$\rightarrow$}\\} \ 
	\xy
	\POS+(0,13)
	\xymatrix"A" @R=0pc @C=1pc @H=1pc @W=0.7pc  
	{
	&& *++[F-]{\phantom{\bullet}}  & \\ 
	&*++[F-]{\phantom{\bullet}} & *++[F-]{\phantom{\bullet}}  & \\ 
	&*++[F-]{\phantom{\bullet}} & *++[F-]{\phantom{\bullet}} &  \\ 
	*++[F-]{\bullet} \ar@{-}[r] & *++[F-]{\bullet} & *++[F-]{\phantom{\bullet}}  & \\ 
	&*++[F-]{\bullet} \ar@{-}[r] & *++[F-]{\bullet}  & 	 \\ 
	& & *++[F-]{\bullet}  & \\ 
	}
	\POS+(40,-26)
	\xymatrix"B"{
	  (3,13)  \\
	 {\bullet} \ar@{-}["A"ddddrr] 
	}
	\POS*\frm{--}
	\endxy
$$


Finally, the composition   $(7,27) A  \To  (5,20) \To (6,24) C $
acts as a left cancellation followed by a right addition.
$$
\xy
		\xymatrix"A" @R=0pc @C=1pc @H=1pc @W=0.7pc  
	{
	& *++[F-]{\phantom{\bullet\ \bullet}}  & \\ 
	*++[F-]{\phantom{\bullet }} & *++[F-]{\phantom{\bullet\ \bullet}}  & \\ 
	*++[F*]{\phantom{\bullet }} & *++[F-]{\phantom{\bullet\ \bullet}}  & \\ 
	*++[F-]{\bullet} & *++[F-]{\bullet\ \bullet}  & \\ 
	}
	\POS+(30,-13)
	\xymatrix"B"{
	  (5,20)  \\
	 {\bullet} \ar@{-}["A"ddr] 
		}
	\POS*\frm{--}
\endxy
	\qquad \xymatrix{\\  \scalebox{3.}{$\rightarrow$}\\} \qquad  
\xy
	\POS+(0,-8)
	\xymatrix"B"{
	  (5,20)  \\
	 {\bullet}  
		}
	\POS*\frm{--}
	\POS+(20,4.5)
	\xymatrix"A" @R=0pc @C=1pc @H=1pc @W=0.7pc  
	{
	 *++[F-]{\phantom{\bullet}}   \\ 
	*++[F-]{\phantom{\bullet }}   \\ 
	*++[F-]{\phantom{\bullet }} \\ 
	*++[F-]{\bullet} \ar@{-}["B"uu] & \\ 
	}
\endxy
$$

\subsection{A preliminary result and a general vision.}  The second 
example of unexpected quotation mentioned in my Bielefeld talk 
deals with  Happel - Ringel's paper  \cite{HR}.  Indeed, the second 
result in the references of  \cite{HR}  has the following form:

\begin{quotation}
 \textbf{[2]} D'Este, G. Talk at Oberwolfach conference on 
representation theory 1981, unpublished.
\end{quotation}

Moreover, the first lines of  \cite[page  164]{HR}  say the following:

\begin{quotation}
``We remark that our account on the decomposition of $\hat C-mod$
  into the module classes   $M_m$ and   $M_{m,m+1}$ 
follows closely the treatment given by  Gabriella d'Este 
in her Oberwolfach talk 1981 [2].'' 
\end{quotation}

I am not able to find the slides with the example described 
in my Oberwolfach talk and/or the handwritten pages 
containing a more general result, suggested by the same 
example.  I remember that I had a long list of hypotheses 
on quivers, algebras and modules into the game, used to 
obtain new algebras with very few new (and not too big)
indecomposable modules.  On the other hand, complexes, 
functors and all the mysteries of derived categories did not 
showed up at all, neither in what I proved nor in what I ever 
hoped to prove.  The best way to sum up all what I can say 
about my Oberwolfach talk 1981  will be to copy, without any 
change, the last slides prepared for my Bielefeld talk 2010 , 
and actually shown at the end.  Almost always I omit the final
part of my talks, but this did not happened on that occasion.  
Indeed, before the beginning of the lectures, the chairman 
opened the files of the various presentations.  When the slide 
of my talk (with a photo of the University of Bielefeld) appeared 
on the screen, Ringel come to my seat.  Since he was leaving 
shortly afterwards, I showed him all my final slides about the 
official/unofficial history mentioned in this note.  Thanks to a 
few pictures and even less words, I needed a very short time 
to illustrate all what I wanted to say.  I was happy to see that 
my recollections were correct.  That's why at the end of my talk 
in Bielefeld I presented the following telegraphic slides.

\section*{What I remember about my talk (Slide -3).}
A result on the support of the  indecomposable
modules over an algebra obtained from a tame 
algebra after two operations: 
\begin{itemize}
\item a  one - point extension  
\item a  one - point coextension  
\end{itemize}
by means of two regular modules (suggested by an
example which gave me the idea that only few new 
indecomposable modules could appear).  

\section*{What I remember BEFORE my talk (Slide -2).}
\begin{itemize}
\item  A long conversation with  Dieter HAPPEL  and 
   Claus RINGEL.
\item Their hint to present only my original example.
\item Their mysterious comments  ?????  (sometimes in 
   German)  that I didn't try to understand. 
\end{itemize}

\section*{What I believe now (Slide -1).}
Behind  my words, calculations done by hands and 
pictures, they were able to  SEE  and/or  IMAGINE  
completely different objects: 
\begin{itemize}
\item quivers with infinitely many vertices 
    (and not just finitely many); 
\item complexes instead of modules;
\item \dotfill
\end{itemize}

\section{Acknowledgments.}    I would like to express my thanks to
``The BIREP  Group'' , to keep the acronym used in the welcome
message to participants contained in the Conference Guide of  
ICRA 2012 .  I am very grateful to all the mathematicians of the 
present and the future generation working and studying in 
Bielefeld for their inspiration and help during international 
conferences, seminars and informal talks.

\end{document}